\documentclass[a4paper,11pt]{article}
\usepackage{amssymb}
\usepackage{amsfonts}
\usepackage{amsmath,tabu}
\usepackage{amsthm}
\usepackage{cite}
\usepackage{amsmath, upgreek}
\usepackage{amssymb}
\usepackage{amscd}
\usepackage{xspace}
\usepackage{verbatim}
\usepackage{xcolor}
	
\usepackage{color}

\newcommand{\s}{s}

\newcommand{\R}{\mathbb{R}}

\newcommand{\N}{\mathbb{N}}

\newcommand{\cuad}{{\sqcap\kern-.68em\sqcup}}

\newcommand{\norm}[1]{\|#1\|}

\numberwithin{equation}{section}
\newtheorem{theorem}{Theorem}[section]

\newtheorem{proposition}[theorem]{Proposition}

\newtheorem{lemma}[theorem]{Lemma}
\newtheorem{corollary}[theorem]{Corollary}
\newtheorem{remark}[theorem]{Remark}
\newcommand{\bremark}{\begin{remark} \em}
\newcommand{\eremark}{\end{remark} }

\newcommand{\cA}{{\mathcal A}}

\newcommand{\cD}{{\mathcal D}}
\newcommand{\cE}{{\mathcal E}}
\newcommand{\cF}{{\mathcal F}}

\newcommand{\cH}{{\mathbb H}}

\newcommand{\cL}{{\mathcal L}}

\newcommand{\cN}{{\mathcal N}}
\newcommand{\cO}{{\mathcal O}}

\newcommand{\lnlap}{L_{\text{\tiny $\Delta \,$}}\!}

\begin{document}
\begin{center}{\bf  \large Bounds for eigenvalues of the Dirichlet problem
\\[2mm]
for the logarithmic Laplacian}\medskip
\bigskip

  {\small  Huyuan Chen\footnote{chenhuyuan@yeah.net}  }
 \medskip

 {\small   Department of Mathematics, Jiangxi Normal University, Nanchang,\\ Jiangxi 330022, PR China   } \\[4mm]

  {\small  Laurent V\'eron\footnote{ veronl@univ-tours.fr }
 \medskip

 {\small   Laboratoire de Math\'{e}matiques et Physique Th\'{e}orique, Universit\'e de Tours, \\ 37200 Tours, France   } \\[3mm]

}

 \medskip

\begin{abstract}
We provide bounds for the sequence of eigenvalues $\{\lambda_i(\Omega)\}_i$ of the Dirichlet problem
$$ \lnlap u=\lambda  u\ \  {\rm in}\ \,  \Omega,\quad\quad  u=0\ \ {\rm in}\  \ \R^N\setminus \Omega,$$
where $\lnlap$ is the logarithmic Laplacian operator  with Fourier transform symbol $2\ln |\zeta|$.
The  logarithmic Laplacian operator is not positively definitive if the volume of the domain is large enough. In this article, we obtain the upper and lower bounds for the sum of the first $k$ eigenvalues by extending the Li-Yau  method and Kr\"oger's method respectively. Moreover,  we show the limit of the sum of the first $k$ eigenvalues, which is independent of the volume of the domain.    Finally, we discuss the lower and upper bounds of the $k$-th principle eigenvalue,  the asymptotic behavior of the limit of eigenvalues. 

 \end{abstract}

\end{center}
 \vspace{1mm}


\vspace{1mm}

\setcounter{equation}{0}
\section{Introduction and main results}

Let $ \lnlap$ be the logarthmic Laplacian in $\R^N$, $N\geq 1$, defined by
\begin{align}
 \lnlap  u(x)&= c_{N} \int_{\R^N  } \frac{ u(x)1_{B_1(x)}(y)-u(y)}{|x-y|^{N} } dy + \rho_N u(x),
    \label{representation-main}
\end{align}
where
\begin{equation}
  \label{eq:def-c-N}
c_N:= \pi^{-  N/2}  \Gamma( N/2) = \frac{2}{\omega_{_{N-1}}}, \qquad \rho_N:=2 \ln 2 + \psi(\tfrac{N}{2}) -\gamma,
\end{equation}
$\omega_{_{N-1}}:=H^{N-1}(S^{N-1})=\int_{_{^{\!S^{N-1}}}}dS$, $\gamma= -\Gamma'(1)$ is the Euler Mascheroni constant and $\psi = \frac{\Gamma'}{\Gamma}$ is the Digamma function. \smallskip

 The aim of this article is to provide estimates of the eigenvalues of the operator $\lnlap$ in a bounded domain $\Omega\subset\R^N$, which are the real numbers $\lambda$ such that there exists a solution to the
 Dirichlet  problem
\begin{equation}\label{eq 1.1}
\left\{ \arraycolsep=1pt
\begin{array}{lll}
 \lnlap u=\lambda  u\quad \  &{\rm in}\quad   \Omega,\\[2mm]
 \phantom{ \lnlap   }
  u=0\quad \ &{\rm{in}}\  \quad \R^N\setminus \Omega.
\end{array}
\right.
\end{equation}

In recent years, there has been a renewed and increasing
interest in the study of boundary value problems involving linear and nonlinear integro-differential operators.
This growing interest is justified both  seminal advances in the
understanding of nonlocal phenomena from a PDE or a probabilistic point of view, see e.g.
\cite{CS0,CS1,CS2,CFQ,EGE,FQ2,GS1,RS,RS1,musina-nazarov} and the references therein, and by important applications.
Among nonlocal differential order operators, the simplest and most studied examples, are the fractional powers of the Laplacian which exhibit many phenomenological  properties.
Recall that for $\s\in(0,1)$ the fractional Laplacian of a function $u \in C^\infty_c(\R^N)$ is defined by
\begin{equation*}
  \label{eq:Fourier-representation}
\mathcal{F}((-\Delta)^\s u)(\xi) = |\xi|^{2\s}\widehat u (\xi)\qquad \text{for all $\xi \in \R^N$},
\end{equation*}
where and in the sequel both $\mathcal{F}$ and $\widehat \cdot$ denote the Fourier transform. Equivalently, $(-\Delta)^\s$ can be written as a singular integral operator under the following form
 \begin{equation}\label{fl 1}
 (-\Delta)^\s  u(x)=c_{N,\s} \lim_{\epsilon\to0^+} \int_{\R^N\setminus B_\epsilon(x) }\frac{u(x)-
u(y)}{|x-y|^{N+2\s}}  dy ,
\end{equation}
where $c_{N,\s}=2^{2\s}\pi^{-\frac N2}\s\frac{\Gamma(\frac{N+2\s}2)}{\Gamma(1-\s)}$ and $\Gamma$  is the  Gamma function, see e.g. \cite{RS1}.

The fractional Laplacian has the following limiting
properties when $\s$ approaches the values $0$ and $1$:
\begin{equation*}
\lim_{\s\to1^-}(-\Delta)^\s u(x)=-\Delta u(x)
\quad \text{and}\quad \lim_{\s\to0^+}  (-\Delta)^\s  u(x) = u(x)\qquad \text{for $u\in C^2_c(\R^N)$,}
\end{equation*}
see e.g. \cite{EGE}.    Recently,    \cite{CT}  shows a further expansion at $s=0$ that for $u \in C^2_c(\R^N)$ and  $x \in \R^N$,
$$
(-\Delta)^\s  u(x) = u(x) + s \lnlap u (x) + o(\s) \quad \text{as\ \, $\s\to0^+$}
$$
 where, formally, the operator
\begin{equation}\label{deriv}
\lnlap:= \frac{d}{d\s}\Big|_{\s=0} (-\Delta)^\s
\end{equation}
 is given as a {\em logarithmic Laplacian}; indeed,
\begin{enumerate}
\item[(i)] for $1 < p \le \infty$, we have $\lnlap  u \in L^p(\R^N)$ and $\frac{(-\Delta)^\s u- u}{\s} \to \lnlap  u$ in $L^p(\R^N)$ as $\s \to 0^+$;
\item[(ii)] $\mathcal{F}(\lnlap u)(\xi) = 2 \ln |\xi|\,\widehat u (\xi)$ 
 \, for a.e. $\xi \in \R^N$.
\end{enumerate}

Note that the problems with integral-differential operators given by kernels with a singularity of order $-N$ have received growing interest recently, as they give rise to interesting limiting regularity properties and Harnack inequalities without scaling invariance, see e.g. \cite{KM}.  Another important domain of study consists in understanding the eigenvalues of the Dirichlet problem with zero exterior value\cite{CT}.  We refer to \cite{JSW,FKT}
for more  topics related to the logarithmic Laplacian and also \cite{JW,FKV13} for general nonlocal operator and related embedding results.  Let $\cH(\Omega)$ denote the space of all measurable functions $u:\R^N\to \R$ with $u \equiv 0$ in $\R^N \setminus \Omega$ and
$$
\int \!\!\!\! \int_{\stackrel{x,y \in \R^N}{\text{\tiny $|x-y|\!\le\! 1$}}} \frac{(u(x)-u(y))^2}{|x-y|^N} dx dy <+\infty.
$$
As we shall see it, $\cH(\Omega)$ is a Hilbert space under the inner product
$$
\mathcal{E}(u,w)=\frac{c_N}2 \int \!\!\int_{\stackrel{x,y \in \R^N}{\text{\tiny $|x-y|\!\le\! 1$}}} \frac{(u(x)-u(y))(w(x)-w(y))}{|x-y|^N}
dx dy,
$$
where $c_N$ is given in (\ref{eq:def-c-N}), with associated norm $\norm{u}_{\cH(\Omega)}=\sqrt{\mathcal{E}(u,u)}$. By \cite[Theorem 2.1]{CP}, the embedding $\cH(\Omega) \hookrightarrow L^2(\Omega)$ is compact. Throughout this article we identify $L^2(\Omega)$ with the space of functions in $L^2(\R^N)$ which vanish a.e. in $\R^N \setminus \Omega$. The quadratic form associated with $\lnlap$ is well-defined on $\cH(\Omega)$ by
$$
\cE_L: \cH(\Omega) \times \cH(\Omega) \to \R, \quad \cE_L(u,w)= \mathcal{E}(u,w) - c_N
\int \!\!\!\! \int_{\stackrel{x,y \in \R^N}{\text{\tiny $|x-y|\!\ge\! 1$}}} \frac{u(x)w(y)}{|x-y|^N} dx dy + \rho_N \int_{\R^N}  u w\,dx,
$$
where $\rho_N$ is defined in (\ref{eq:def-c-N}). A function $u \in \cH(\Omega)$ is an eigenfunction of (\ref{eq 1.1}) corresponding to the eigenvalue
$\lambda$ if
$$
\cE_L(u,\phi) = \lambda \int_{\Omega}u\phi\,dx \qquad \text{for all $\phi \in \cH(\Omega)$.}
$$

\begin{proposition}\label{teo 1}\cite[Theorem 1.4]{CT}
Let $\Omega$ be a bounded domain  in $\R^N$.  Then problem (\ref{eq 1.1}) admits a sequence of eigenvalues
$$
\lambda_1 (\Omega)<\lambda_2 (\Omega)\le \cdots\le \lambda_i (\Omega)\le \lambda_{i+1} (\Omega)\le \cdots
$$
and corresponding eigenfunctions $\phi_i$, $i \in \N$ such that the following holds:
\begin{enumerate}
\item[(a)] $\lambda_{i}(\Omega)=\min \{\cE_L(u,u) \::\:  u\in \cH_i(\Omega)\::\: \norm{u}_{L^2(\Omega)}=1\}$, where
$$
\cH_1(\Omega):= \cH(\Omega)\quad \text{and}\quad  \cH_i(\Omega):=\{u\in\cH(\Omega)\::\: \text{$\int_{\Omega} u \phi_i \,dx =0$ for $i=1,\dots i-1$}\}\quad \text{for $i>1$;}
$$
\item[(b)] $\{\phi_i\::\: i \in \N\}$ is an orthonormal basis of $L^2(\Omega)$;
\item[(c)] $\phi_1$ is positive in $\Omega$. Moreover, $\lambda_1 (\Omega)$ is simple, i.e., if $u \in \cH(\Omega)$ satisfies (\ref{eq 1.1}) in weak sense with $\lambda = \lambda_1 (\Omega)$, then $u=t\phi_1$ for some $t\in\R$;
\item[(d)] $\lim \limits_{i\to \infty} \lambda_i (\Omega)=+\infty$.
\end{enumerate}
\end{proposition}

Due to lack of the homogenous property for   the logarithmic Laplacian operator,   the effect of  the domain for the principle eigenvalue can't be expected as the Laplacian or fractional Laplacian, just by scaling the domain by their homogeneous property of such operators. Secondly, the logarithmic Laplacian operator is no longer positively definitive if $|\Omega|$ is to large, since it is proved in \cite{CT} that the positivity of the principle eigenvalue is equivalent to the comparison principle, which  does not hold for  balls with large radius.  
 These properties of  the logarithmic Laplacian operator  enrich the asymptotics of the Dirichlet eigenvalues as we will see below, but also make more difficult the obtention of bounds for eigenvalues.  

It is well-known that the Hilbert-P\'olya conjecture is to associate the zero of the Riemann Zeta function with the eigenvalue of a Hermitian operator. This quest initiated the mathematical interest  for estimating the sum of Dirichlet eigenvalues of the Laplacian while in physics the question is related to count the number of bound states of a one body Schr\"odinger operator and to study their asymptotic distribution. In 1912, Weyl in \cite{W} shows that  the $k$-th eigenvalue $\mu_k(\Omega)$ of Dirichlet problem with the Laplacian operator
\begin{equation}\label{eq 1.1-lap}
\left\{ \arraycolsep=1pt
\begin{array}{lll}
-\Delta u=\mu  u\quad \  &{\rm in}\quad   \Omega,\\[2mm]
 \phantom{ \lnlap   }
  u=0\quad \ &{\rm{in}}\  \quad \R^N\setminus \Omega
\end{array}
\right.
\end{equation}
 has the asymptotic behavior
$\mu_k(\Omega)\sim C_N(k|\Omega|)^{\frac2N}$ as $k\to+\infty$, where $C_N=(2\pi)^2|B_1|^{-\frac{2}{N}}$. Later, P\'olya \cite{P} (in 1960) proved that
\begin{equation}\label{p-conj}
\mu_k(\Omega)\geq C\left(\frac{k}{|\Omega|}\right)^{\frac{2}N}
\end{equation}
holds for $C=C_N$ and any "plane-covering domain" $D$ in $\R^2$, (his proof also works in dimension $N\geq 3$)
 and he also conjectured that (\ref{p-conj}) holds with $C=C_N$ for any bounded domain in $\R^N$. Rozenbljium \cite{Ro} and independently Lieb \cite{L} proved (\ref{p-conj}) with a positive constant $C$ for general bounded domain.   Li-Yau \cite{LY} improved
 the constant $C=\frac{N}{N+2}C_N$, and with that constant (\ref{p-conj}) is also called Berezin-Li-Yau inequality because this constant is achieved with the help of Legendre transform as in the Berezin's earlier paper \cite{Be}.  The Berezin-Li-Yau inequality then is generalized in \cite{CP,CgW,L,K,M,CgY},   for degenerate elliptic operators in \cite{CZ,HY,YY} for the fractional Laplacian $(-\Delta)^s$ defined in (\ref{fl 1}) and the inequality reads
 \begin{equation}\label{p-conjf}
\mu_{s,k}(\Omega)\geq \frac{N}{N+2s}C_N\left(\frac{k}{|\Omega|}\right)^{\frac{2s}N}.
\end{equation}
 Due to the expression of the  Fourier symbol of $\lnlap$,  Berezin-Li-Yau method can not be applied to our problem (\ref{eq 1.1}). Our results  are based on the appropriate  estimates for the solutions of equations:
 $$r\ln r =c\quad\ {\rm and}\quad\ \frac{r}{\ln r-\ln\ln r}=t.$$
The  estimates that we obtain  provide a uniform lower bound  of the sum of  the first $k$-eigenvalues, independently of $k$, an estimate which has a particular interest when these eigenvalues are negative. More precisely, we have the following inequalities:

\begin{theorem}\label{teo 1.3}
Let $\Omega\subset \R^N$be a bounded domain, $\{\lambda_i(\Omega)\}_{i\in\N}$ be the sequence of eigenvalues of  problem (\ref{eq 1.1}) obtained in Proposition \ref{teo 1} and define
\begin{equation}\label{dn-1}
d_N=\frac{2\omega_{_{N-1}}}{N^2(2\pi)^N}
\end{equation}
   Then there holds \smallskip

\noindent (i) for any $k\in\N^*$,
$$
 \sum^k_{i=1}\lambda_i(\Omega)\geq -d_N |\Omega|;$$
 $(ii)$ if $k>\frac{eN d_N}2  |\Omega| $,
 $$\sum^k_{i=1}\lambda_i(\Omega)>0;$$
 $(iii)$   if $k\geq    \frac{e^{e+1}Nd_N}{2}|\Omega|$,
 \begin{equation}\label{lower bb}
 \sum^k_{i=1}\lambda_i(\Omega) \geq
 \frac{2k}{N} \left(\ln k+  \ln\Big(\frac{2}{eNd_N|\Omega|}\Big) -\ln\ln\Big(\frac{2k}{eNd_N|\Omega|}\Big)\right).
 \end{equation}

 \end{theorem}

Our second interest is to give an upper bound for the first $k$ eigenvalues. 
Motivated by Kr\"oger's result  for the Laplacian  \cite{K}, we shall build an upper bound by
 calculating the related Rayleigh quotient via a particular complex valued function.  Together with the lower bound (\ref{lower bb}), we can derive the limit of the sum of the first $k$ eigenvalues  as $k\to+\infty$. The results state as following.


\begin{theorem}\label{teo 1.1-sum}
Let $\Omega\subset \R^N$ be a bounded domain and   $\{\lambda_{i}(\Omega)\}_{i\in\N}$ be the sequence of eigenvalues of  problem (\ref{eq 1.1}).
Then  for $k> \frac{e Nd_N}{2}|\Omega| $,
  \begin{equation}\label{upper bb}
  \sum^k_{i=1}\lambda_{i}(\Omega) \leq  \frac{2k}{N} \left(\ln (k +1)+\ln \Big(\frac{p_N}{|\Omega|}\Big)+ \frac{\omega_{_{N-1}}}{\sqrt{|\Omega|}}  \ln \ln \Big(\frac{p_N (k+1)}{|\Omega|}\Big) \right)
  \end{equation}
  and
  \begin{equation}\label{upper bb-lim}
 \lim_{k\to+\infty} (k\ln k)^{-1} \sum^k_{i=1}\lambda_{i}(\Omega) =  \frac{2}{N},
  \end{equation}
where $p_N=\frac{2 (2\pi)^NN}{\omega_{_{N-1}} }$.

 \end{theorem}

Note that the assumption that $k> \frac{e Nd_N}{2}|\Omega|$ is required to make sure that 
$\lambda_{k_0}>0$, here $k_0$ is the smallest  positive integer $k_0\geq \frac{e Nd_N}{2}|\Omega|$.

From the bounds of sum of eigenvalues from our main results,  we can provide the following  the Wely's formula for the  logarithmic Laplacian. 


 \begin{corollary}\label{cr 1.2}
Let $\Omega\subset \R^N$ be a bounded domain and   $\{\lambda_{i}(\Omega)\}_{i\in\N}$ be the sequence of eigenvalues of  problem (\ref{eq 1.1}).
Then 
$(i)$
 \begin{equation}\label{Weyl-limit}
  \lim_{k\to+\infty}\frac{\lambda_{k}(\Omega)}{ \ln k } =\frac2N.
  \end{equation}
  $(ii)$
  \begin{equation}\label{Weyl-upplow}
  \frac{2}{N} \Big(\ln k+\ln \frac{2}{eN d_N|\Omega|}\Big)\leq  \lambda_k(\Omega)\leq  \frac{2}{N} \ln k+c_0\big(\ln\ln (k+e)\big)^2+2\ln \frac{|B_1|}{|\Omega|}.
    \end{equation}
    where $c_0>0$ is independent of $k$ and $\Omega$.
 \end{corollary}

\begin{remark}
$(a)$ The limit of $\frac{\lambda_{k}(\Omega)}{ \ln k }$  and the one of $(k\ln k)^{-1} \sum^k_{i=1}\lambda_{i}(\Omega)$ as $k\to+\infty$ have the same value $\frac{2}{N}$, which is independent of $\Omega$;\smallskip

 $(b)$   The Weyl formula (\ref{Weyl-limit}) also could be proved by 
$$\lim_{\lambda\to+\infty} e^{-\frac{\lambda}2N}\cN(\lambda) =\frac{|\Omega|\omega_{_{N-1}}}{N(2\pi)^N}, $$
which is shown  \cite[Corollary 6.2]{LW}, where $\cN$ is the counting function of eigenvalues, 
$$
\cN(t) = \sum_{j \in \N}sgn_+( t - \lambda_j(\Omega) )=\sum_{j \in \N}( t - \lambda_j(\Omega) )_+^0.
$$
Here $sgn_+(r)=1$ if $r>0$, $sgn_+(r)=0$ if $r\leq 0$ and $r_\pm = (|r|\pm r)/2$ denote the positive and negative part of  $x \in \R$. More estimates for $\cN$ see \cite{LW}.

$(c)$ The limit (\ref{Weyl-limit})  is proved  by showing the inequalities  
 \begin{equation}\label{lower bb-1}
  \ln k+ \ln \frac{2}{eNd_N|\Omega|} -\ln\ln \frac{2k}{eNd_N|\Omega|}  \leq  \frac{N}2 \lambda_k(\Omega)\leq   \ln k+2 + \tilde c_0\big(\ln\ln (k+1)\big)^2
  \end{equation}
 for $k\geq  \frac{e^{e+1}Nd_N}{2}|\Omega|$, where $\tilde c_0>0$ is independent of $k$, but dependent of $|\Omega|$. Here the first and second inequalities  follow by (\ref{lower bb}) and (\ref{upper bb}) respectively, along with the monotonicity of the sequence of eigenvalues. 
  
 The inequalities  (\ref{Weyl-upplow}) provide sharper bounds for $\lambda_k(\Omega)$, by the aid of a scaling property and some estimates of $\cN$ in \cite{LW};
\smallskip

 \end{remark}


  The rest of this paper is organized as follows.  Section 2 is devoted to proving the lower bound  by developing Li-Yau's method, and then we prove Theorem \ref{teo 1.3}. In Section 3, we show the upper bounds for the first $k$-eigenvalues in Theorem \ref{teo 1.1-sum}.   Finally, we prove the Wely's limit of eigenvalues in Corollary \ref{cr 1.2} and  obtain the more bounds for $\lambda_k(\Omega)$ in Section 4.


\section{Lower bounds}

Let
$$g(r)=r\ln r\quad{\rm for}\ \, r>0,$$
then
$g(e)=e$, $g(1)=0$ and $g(\frac1e)=-\frac1e$.

 \begin{lemma}\label{lm 2.1}
  For $c\geq -\frac1e$, there exists  a unique point $r_c\geq \frac1e$  such that
$$g(r_c)=c,$$
and we have that $ r_c\leq 1+c.$ Furthermore,  \smallskip

 \noindent $(i)$ for $-\frac1e\leq c\leq 0$,
 $$r_c\geq  1+(e-1)c \geq \frac1e;$$
 \noindent $(ii)$ for $0\leq c\leq e$,
 $$r_c\geq 1+\frac{e-1}{e}c;$$
\noindent $(iii)$ for $c\geq e$,
$$r_c\geq 1+\frac{e-1}{e}c
$$
and
 \begin{equation}\label{N1}
 \frac{c}{\ln c } \leq r_c\leq \frac{c}{\ln c- \ln \ln c}.
 \end{equation}

 \end{lemma}
 \noindent{\bf Proof.} The function $g$ is increasing in $[\frac1e,+\infty)$ with value in $[-\frac1e,+\infty)$. Hence $r_c$ is uniquely determined if $c\geq -\frac1e$, $c\mapsto r_c$ is increasing from $[-\frac1e,+\infty)$ onto $[\frac1e,+\infty)$,
and $g$ is convex.\smallskip

For $a>0$, we define $\psi_a(x)=(1+ax)\ln(1+ax)-x$ for $x>-\frac 1a$. Then $\psi_a(x)>0$ (resp. $\psi_a(x)<0$) is equivalent to $1+ax>r_x$ (resp. $1+ax<r_x$).
Note that $\psi'_a(x)=a(1+\ln(1+ax))-1$. Since $\psi'_a(-\frac 1a)=-\infty$ and $\psi'_a$ is increasing,
 $\psi'_a(0)=a-1$ is the maximal (resp. minimal) value of $\psi'_a$ on $(-\frac 1a,0]$ (resp. on $[0,\infty)$). Therefore, if $a>1$,
$\psi_a$ is positive on $(-\frac{1}{a},r^*_a)$ for some $r^*_a\in (-\frac{1}{a},0)$, negative on $(r^*_a,0)$ and positive on $(0,\infty)$. If $0<a<1$, $\psi_a$ is positive on $(-\frac{1}{a},0)$, negative on $(0,r^*_a)$ for some $r^*_a>0$ and positive on $(r^*_a,\infty)$. If $a=1$, $\psi_1$ is positive on $[-\frac 1e,0)\cup (0,\infty)$ and vanishes only at $0$.
Then $\psi_1\geq 0$ implies the first assertion.

Since $e-1>1$ and $\psi_{e-1}(-\frac{1}{e})=0$, $\psi_{e-1}(x)<0$ for $x\in (-\frac{1}{e},0)$. This gives (i).

Since $0<\frac{e-1}{e}<1$, $\psi_{\frac{e-1}{e}}$ is negative on $(0,r^*_\frac{e-1}{e})$ and positive on $(r^*_\frac{e-1}{e},\infty)$. Since
$\psi_{\frac{e-1}{e}}(e)=0$, $r^*_\frac{e-1}{e}=e$ and we get (ii) and (iii).

 Since $g$ is increasing on $[e,\infty)$, \eqref{N1} is equivalent to
$$c-\frac{\ln c}{\ln\ln c}\leq c\leq \frac{c}{\ln c-\ln\ln c}\ln\left(\frac{c}{\ln c-\ln\ln c}\right)
=c\frac{\ln c-\ln(\ln c-\ln\ln c)}{\ln c-\ln\ln c}.
$$
Set $C=\ln c$, then
$$\frac{\ln c-\ln(\ln c-\ln\ln c)}{\ln c-\ln\ln c}=\frac{C-\ln(C-\ln C)}{C-\ln C}>1\quad\text{for }\; C>1
$$
and \eqref{N1} follows.\hfill$\Box$ \medskip

  \begin{lemma}\label{lm 2.2}
  Let $f$ be a real-valued function defined in $\R^N$ with $0\leq f\leq M_1$ and
  $$2\int_{\R^N} \ln|z|\, f(z)dx =M_2.$$
 Then $(i)$
 $$M_2\geq -\frac{2\omega_{_{N-1}}}{N^2}M_1;$$
  $(ii)$
  $$\int_{\R^N} f(z)dz \leq  \frac{M_1\omega_{_{N-1}}}{N}\left(e+\frac{N^2M_2}{2M_1\omega_{_{N-1}}}\right)=\frac{e\omega_{_{N-1}}}{N}M_1+\frac N2M_2; $$
   $(iii)$   assuming more that $\frac{M_2}{M_1}\geq \frac{2e^2\omega_{_{N-1}}}{N^2}$, there holds
  $$\int_{\R^N}f(z)dz\leq \frac{NM_2}{2}\left(\ln\left(\frac{N^2M_2}{2eM_1\omega_{_{N-1}}}\right)-\ln\ln\left(\frac{N^2M_2}{2eM_1\omega_{_{N-1}}}\right)\right)^{-1}. $$

  \end{lemma}
 \noindent{\bf Proof. }
We have
$$\begin{array}{lll}\displaystyle
\frac{M_2}{2}=\int_{B_1}\ln|z| f(z)dz+\int_{B^c_1}\ln|z| f(z)dz\\[4mm]
\phantom{\frac{M_2}{2}}\displaystyle \geq M_1\int_{B_1}\ln|z| dz+\int_{B^c_1}\ln|z| f(z)dz
\ \geq -\frac{\omega_{_{N-1}}}{N^2}M_1.
\end{array}$$
Hence (i) holds.

For $R>0$ we have that
 $$(\ln|z|-\ln R)(f(z)-M_1{\bf 1}_{B_R})\geq 0.
 $$
By integration over $\R^N$ we get
 $$\frac{M_2}{2}+\frac{M_1\omega_{_{N-1}}R^N}{N^2}\geq \ln R\int_{\R^N}f(z)dz.
 $$
The estimate from above of $\int_{\R^N}f(z)dz$ is obtained by
 \begin{equation}\label{O-1}
\int_{\R^N}f(z)dz\leq\inf\Big\{A>0 \text{ s.t. }\frac{M_2}{2}+\frac{M_1\omega_{_{N-1}}R^N}{N^2}-A\ln R\geq 0\text{ for all }R>0\Big\}.
\end{equation}
Set
$$\Theta_A(R)=\frac{M_2}{2}+\frac{M_1\omega_{_{N-1}}R^N}{N^2}-A\ln R,
$$
then $\Theta_A$ achieves the minimum if
$$\frac{M_1\omega_{_{N-1}}R^N}{N}=A\Longleftrightarrow R=R_A:=\left(\frac{NA}{\omega_{_{N-1}}M_1}\right)^{\frac{1}{N}}.$$
Hence
 \begin{equation}\label{O-2}
\Theta_A(R_A)=\frac{M_2}{2}+\frac{A}{N}-\frac{A}{N}\ln\left(\frac{NA}{\omega_{_{N-1}}M_1}\right).
\end{equation}
Put $r=\frac{NA}{M_1\omega_{_{N-1}}}$, then
 \begin{equation}\label{O-3}
 \Theta_A(R_A)\geq 0\Longleftrightarrow r\ln r-r\leq \frac{N^2M_2}{2M_1\omega_{_{N-1}}}
 \Longleftrightarrow g\left(\frac re\right)\leq \frac{N^2M_2}{2eM_1\omega_{_{N-1}}}.
\end{equation}
Then  $\tfrac{r}{e}\leq r_c$ with $c=\frac{N^2M_2}{2eM_1\omega_{_{N-1}}}$, inequality $r_c\leq 1+c$ in Lemma \ref{lm 2.1} yields
$$r=\frac{NA}{M_1\omega_{_{N-1}}}\leq e+\frac{N^2M_2}{2M_1\omega_{_{N-1}}}\Longrightarrow
\int_{\R^N}f(z)dz\leq \frac{M_1\omega_{_{N-1}}}{N}\left(e+\frac{N^2M_2}{2M_1\omega_{_{N-1}}}\right),
$$
which is (ii).
\smallskip

Assuming now that $\frac{M_2}{M_1}\geq \frac{2e^2\omega_{_{N-1}}}{N^2}$, we can apply Lemma \ref{lm 2.1}-(iii) and get
$$\int_{\R^N}f(z)dz\leq \frac{NM_2}{2}\left(\ln\left(\frac{N^2M_2}{2eM_1\omega_{_{N-1}}}\right)-\ln\ln\left(\frac{N^2M_2}{2eM_1\omega_{_{N-1}}}\right)\right)^{-1},
$$
which is (iii) and ends the proof.
\hfill$\Box$ \medskip

 \begin{lemma}\label{lm 2.3}
 Let
$$\tilde g(r)=\frac{r}{\ln r-\ln\ln r }  \quad{\rm for}\ \, r>e.$$
Then for $t>\frac{e^e}{e-1}$, there exists  a unique point $r_t>e$ such that $\tilde g(r_t)=t.$
Furthermore,
 \begin{equation}\label{O-4}
 t(\ln t-\ln\ln t)\leq r_t <   t\ln t. \end{equation}
 \end{lemma}

\noindent{\bf Proof}. Since
\begin{align*}
\tilde g' (r)&=\frac{1}{\ln r-\ln\ln r} -\frac{1-  (\ln r)^{-1}}{(\ln r-\ln\ln r )^2} \\
&\geq \frac{1}{\ln r-\ln\ln r }\Big(1-\frac{1}{\ln r-\ln\ln r}\Big)>0,
\end{align*}
the function $\tilde g$ is increasing from $(e,+\infty)$ onto $(\frac{e^e}{e-1},+\infty)$. Setting $r^*_t=t(\ln t-\ln\ln t)$,  then
 \begin{align*}
\tilde g(r^*_t)&= \frac{t(\ln t-\ln\ln t)}{\ln t+\ln(\ln t-\ln\ln t) -\ln\ln(t(\ln t-\ln\ln t)) }  \\
&\leq \frac{t(\ln t-\ln\ln t)}{\ln t+\ln(\ln t-\ln\ln t) -\ln\ln(t\ln t ) } \\
&=\frac{  \ln t-\ln\ln t}{\ln t  -\ln\frac{\ln(t\ln t )}{\ln t-\ln\ln t} }t
\\&\leq  t,
\end{align*}
where the last inequality holds if
 \begin{align*}
 \frac{\ln(t\ln t )}{\ln t-\ln\ln t} \leq \ln t,
\end{align*}
which is equivalent to
$$\tilde h(\tau):=\tau^2-(\ln\tau+1)\tau-\ln\tau\geq0,\quad \tau=\ln t.$$
Freezing the coefficient $\ln\tau$, $\tilde h(\tau)=(\tau-\tau_1)(\tau-\tau_2)$, where the $\tau_1,\, \tau_2$ depend of $\tau$, but  $\tau_1<0<\tau_2$, since
$\tau_1\tau_2=-\ln\tau<0$. Because $\tilde h(e)=e^2-2e-1=0.9584\pm 10^{-4}$, we have $e>\tau_1$. Hence $\tau>e$ implies $\tau>\tau_1$ which in turn implies
$\tilde h(\tau)>0$. Hence $r^*_t\leq r_t$ using the monotonicity of $\tilde g$.
\smallskip

 Let $s_t=  t\ln t$, then
 \begin{align*}
\tilde g(s_t)&= \frac{t\ln t }{\ln t+\ln\ln t  -\ln\ln (t\ln t) }  < t
\end{align*}
by the fact that
$$\ln\ln t  -\ln\ln (t\ln t)<0\quad{\rm for}\  \, t>e.$$
  Hence $s_t \geq r_t$, which ends the proof.   \hfill  $\Box$

  \medskip

\noindent {\bf Proof of Theorem \ref{teo 1.3}.}
 Denote
 $$\Phi_k(x,y)=\sum^k_{j=1} \phi_j(x)\phi_j(y),\qquad(x,y)\in\Omega\times \Omega,$$
 and
 $$\widehat{\Phi}_k(z,y)=(2\pi)^{-\frac N2} \int_{\R^N} \Phi_k(x,y)e^{ix\cdot z}dx,$$
 where $\widehat{\Phi}_k$ is the Fourier transform with respect to $x$.
 Hence we have that
 \begin{align*}
 \int_{ \R^N}  \int_{\Omega} |\widehat{\Phi}_k(z,y)|^2dzdy =  \int_{\Omega} \int_{\Omega}|\Phi_k(x,y)|^2dxdy=k
 \end{align*}
 by the orthonormality of  the $\{\phi_j\}_{i\in\N}$ in $L^2(\Omega)$. Furthermore, we note that
\begin{equation}\label{B1}\begin{array}{lll}
\displaystyle \int_{\Omega} |\widehat{\Phi}_k(z,y)|^2dy=\int_{\Omega} \left(\sum^k_{j=1}\widehat\phi_j(z)\phi_j(y)\right)\left(\sum^k_{j=1}\overline{\widehat\phi_j(z)}\phi_j(y)\right)dy\\[4mm]
 \phantom{ \int_{\Omega} |\widehat{\Phi}_k(z,y)|^2dy}\displaystyle
 =\int_{\Omega} \left(\sum_{j,\ell=1}^k\widehat\phi_j(z)\overline{\widehat\phi_\ell (z)}\phi_j(y)\phi_\ell(y)\right)dy\\[4mm]
 \phantom{ \int_{\Omega} |\widehat{\Phi}_k(z,y)|^2dy}\displaystyle
 =\sum^k_{j=1}|\widehat\phi_j(z)|^2.
 \end{array}\end{equation}
Using again the orthonormality of  the $\{\phi_j\}_{i\in\N}$ in $L^2(\Omega)$,  we infer by the k-dim Pythagore theorem,
\begin{equation}\label{B2}\begin{array}{lll}
 \displaystyle \int_{\Omega} |\widehat{\Phi}_k(z,y)|^2dy=(2\pi)^{-N}\int_{\Omega}\left|\sum_{j=1}^k\left(\int_{\Omega}e^{ix.z}\phi_j(x)dx\right)\phi_j(y)\right|^2dy\\[4mm]
 \phantom{\int_{\Omega} |\widehat{\Phi}_k(z,y)|^2dy} \displaystyle
 =(2\pi)^{-N}\sum_{j=1}^k\left|\int_{\Omega}e^{ix.z}\phi_j(x)dx\right|^2\\[4mm]
 \phantom{\int_{\Omega} |\widehat{\Phi}_k(z,y)|^2dy} \displaystyle
 \leq (2\pi)^{-N}|\Omega|.
 \end{array}\end{equation}

We have,  from the Fourier expression of $\lnlap$,
 \begin{align*}
\sum^k_{j=1}\lambda_j(\Omega)&=\int_{\Omega}\int_\Omega \Phi_k(x,y) \lnlap\Phi_k(x,y) dy dx\\
&=2\sum^k_{j=1}\int_{\R^N}|\widehat\phi_j(z)|^2\ln|z|dz\\
&=2\int_{\R^N}\left(\int_\Omega|\widehat\Phi_k(z,y)|^2dy\right)\ln|z|dz.
 \end{align*}
 Now we apply Lemma \ref{lm 2.2} to the function
 $$f(z)=\int_{\Omega} |\widehat{\Phi}_k(z,y)|^2dy$$
 with
 $$M_1=(2\pi)^{-N}|\Omega|\quad{\rm and}\quad M_2=\sum^k_{j=1}\lambda_j(\Omega).$$
\\[1mm]
Part $(i)$: By  Lemma \ref{lm 2.2} $(i)$,
 $$\sum^k_{j=1}\lambda_j(\Omega)\geq -\frac{2\omega_{_{N-1}}}{N^2(2\pi)^{N}}|\Omega|=-d_N|\Omega|,$$
 where  $d_N$ is the constant defined in (\ref{dn-1}).\\[1mm]
Part $(ii)$:
\begin{align*}
 k= \int_{\R^N} f(z) dz\leq \frac{e\omega_{_{N-1}}|\Omega|}{N(2\pi)^N} +\frac N2\sum^k_{j=1}\lambda_j(\Omega),
 \end{align*}
 which implies that
 $$\sum^k_{j=1}\lambda_j(\Omega)\geq \frac {2k}{N}-\frac{2e\omega_{_{N-1}}|\Omega|}{N^2(2\pi)^N}.$$
 Part $(iii)$: for $k\in\N$, if
 $$\sum^k_{j=1}\lambda_j(\Omega)\geq \frac{2e^2\omega_{_{N-1}}}{N^2}|\Omega|,$$
 then
  \begin{align*}
  k\leq  \frac{N M_2 }{2  }\left(\ln \left(\frac{N^2M_2 }{2e M_1 \omega_{_{N-1}}}\right)- \ln\ln \left(\frac{N^2M_2}{2e M_1 \omega_{_{N-1}}} \right)\right)^{-1}.
 \end{align*}
Setting
$$r=\frac{N^2M_2}{2e M_1 \omega_{_{N-1}}}\quad{\rm and}\quad t=\frac{Nk}{eM_1\omega_{_{N-1}}}=\frac{(2\pi)^NNk}{e\omega_{_{N-1}}|\Omega|},$$
we have  from \eqref{O-4}  that
\begin{equation}\label{4.1}
r\geq r_t\geq t\left(\ln t-\ln\ln t\right)
\end{equation}
for any $t>e^e$, i.e.
 $$k>\frac{e^{e+1}\omega_{_{N-1}}|\Omega|}{(2\pi)^NN}=\frac{e^{e+1}Nd_N}{2}|\Omega|.$$
This implies
 $$\frac{N^2M_2}{2e M_1 \omega_{_{N-1}}}\geq \frac{(2\pi)^NNk}{e\omega_{_{N-1}}|\Omega|}\left(\ln\Big(\frac{(2\pi)^NNk}{e\omega_{_{N-1}}|\Omega|}\Big)-\ln\ln\Big(\frac{(2\pi)^NNk}{e\omega_{_{N-1}}|\Omega|}\Big)\right),$$
from what we infer
\begin{equation}\label{4.x} \sum^k_{j=1}\lambda_j(\Omega) \geq \frac{2k}{N} \left(\ln\Big(\frac{2k}{eNd_N|\Omega|}\Big)-\ln\ln\Big(\frac{2k}{eNd_N|\Omega|}\Big)\right),
\end{equation}
 which completes the proof.\hfill$\Box$\medskip

\begin{remark} 
 
By the nondecreasing monotonicity of $k\mapsto \lambda_k(\Omega)$, we have that 
for  $k\geq \frac{e Nd_N}{2}|\Omega|$
 $$\lambda_k(\Omega)\geq  \frac{1}{k}\sum^k_{i=1}\lambda_i(\Omega)>0.$$ 
\end{remark}

 
\section{Upper bound }

For any bounded complex valued functions $u,v$ defined on $\Omega$, there holds
\begin{equation}\label{2.0-log-sum}
\lnlap (uv) (x) =u(x) \lnlap v (x) +c_N   \int_{B_1(x)} \frac{u(x)-u(\zeta) }{|x-\zeta|^{N }} v(\zeta)  d\zeta.
\end{equation}

\begin{lemma}\label{lm1-log-sup}
For $z\in \R^N$, we denote
$$\mu_z(x)=e^{{\rm i} x\cdot z},\quad\forall\, x\in\R^N,$$
then
\begin{equation}\label{2.1-log-sum}
\lnlap \mu_z(x)=  (2\ln |z| ) \mu_z(x),\quad\forall\, x\in\R^N.
\end{equation}
\end{lemma}

\noindent{\bf Proof. } {\it Step 1:
we claim that for}
\begin{equation}\label{2.1A}
(-\Delta)^s \mu_z(x)=   |z|^{2s} \mu_z(x),\quad\forall\, x\in\R^N.
\end{equation}

 Without loss of generality, it is enough to
prove (\ref{2.1A}) with $z=te_1$, where $t>0$ and $e_1=(1,0,\cdots,0)\in\R^N$. For this,
we  write
$$v_t(x) =\mu_z(x_1)=e^{{\rm i} tx_1 },\quad x=(x_1,x')\in \R\times \R^{N-1}.$$
 For $N\geq 2$ it implies by \cite[Lemma 3.1]{CV} that
  \begin{eqnarray*}
(-\Delta)^s v_t(x)=(-\Delta)^s_{\R} v_t(x_1).
\end{eqnarray*}
 Now we claim  that
\begin{equation}\label{2.2}
(-\Delta)^s_{\R} v_t(x_1)=t^{2s} v_t(x_1),\quad\forall\, x_1\in\R.
\end{equation}
 Indeed, observe that $-\Delta_{\R}:=- (v_t)_{x_1x_1} =t^2 v_t$ in $\R$ and then
$$  (|\xi_1|^{2}-t^2) \widehat{v_t}(\xi_1)=\cF\left(-\Delta_{\R} v_t-t^2v_t\right)(\xi_1)=0, $$
which implies that
$${\rm supp}(\widehat{v_t})\subset \{\pm t\},$$
which in turn implies
$$(|\xi_1|^{2s}-t^{2s}) \widehat{v_t}(\xi_1)=0=\cF\left((-\Delta)^s_{\R} v_t-t^{2s}v_t\right)(\xi_1). $$
and finally
$$\left((-\Delta)^s_{\R} v_t-t^{2s}v_t\right)(\xi_1)=0\quad{\rm in}\ \ \R,$$
which yields
\begin{align*}
(-\Delta)^s v_t(x) &=  (-\Delta)^s_{\R} v_t =  t^{2s}v_t(x),\quad \forall\,x\in\R^N,
 \end{align*}

\noindent{\it Step 2: we show   (\ref{2.1-log-sum}). } From the property \eqref{deriv} of $\lnlap$ since $\mu_z$ is bounded,
\begin{align*}
 0 &= \frac{(-\Delta)^s\mu_z(x)-|z|^{2s}\mu_z(x) }{s}
 \\&= \frac{(-\Delta)^s\mu_z(x)- \mu_z(x) }{s}- \frac{ |z|^{2s} - 1 }{s}\mu_z(x)
  \\& \to \lnlap \mu_z(x) -(2\ln|z|) \mu_z(x)\qquad  {\rm as}\ \   s\to0^+,
\end{align*}
hence,
$$\lnlap \mu_z(x)=(2\ln|z|) \mu_z(x),\quad\forall\, x\in\R^N,$$
which is the claim.  \hfill$\Box$\medskip

 Next, let $\eta_0\in C^1(\R)$ be a nondecreasing real value function such that $\norm{\eta'_0}_{L^{\infty}}\leq 2$ satisfying
 $$\eta_0(t)=1\ \ {\rm if} \ \, t\geq 1,\qquad \eta_0(t)=0\ \ {\rm if} \ \, t\leq 0.$$
Since $\Omega$ is a bounded domain, there exists a $C^1$ domain $\cO\subset\Omega$ such that $|\cO|\geq  \frac34 |\Omega|$.
 For $\sigma>0$, we set again
\begin{equation}\label{ws-1-sum}
w_\sigma(x)=\eta_0(\sigma^{-1} \bar\rho(x)),\quad \forall\,  x\in\R^N.
\end{equation}
where $\bar\rho(x)=  {\rm dist}(x,\partial\cO)$.
Observe that  $w_\sigma\in \cH_0(\Omega)$ and
$$w_\sigma \to 1 \quad{\rm in}\  \,\cO \ \ {\rm as}\ \ \sigma\to0^+.$$
Thus,   there exists $\sigma_1>0$ such that for $\sigma\in(0,\sigma_1]$,
$$|\Omega|>\int_{\Omega} w_\sigma\, dx\geq \int_{\Omega} w_\sigma^2\, dx>\frac{|\Omega|}2.$$


\begin{lemma}\label{lm 2.2x} Let
 $$\cL_z w_\sigma(x)=\int_{B_1(x)} \frac{w_\sigma(x)-w_\sigma(\zeta)  }{|x-\zeta|^{N }} e^{i\zeta\cdot z} d\zeta, $$
 then  there holds
$$\Big|\cL_z w_\sigma(x) \Big| \leq  \frac{2 \omega_{_{N-1}}}{\sigma}  \quad{\rm for}\ \, x\in \Omega.$$
\end{lemma}

\noindent{\bf Proof.}
Actually, if $x\in \Omega$, we have that
\begin{align*}
| w_\sigma(x)-w_\sigma(\zeta) | \leq  \|Dw_\sigma\|_{L^\infty} |x-\zeta|
 \leq  \sigma^{-1}\|\eta'_0\|_{L^\infty} |x-\zeta|,
\end{align*}
then
\begin{align*}
  \Big|\int_{B_1(x)} \frac{w_\sigma(x)-w_\sigma(\zeta)  }{|x-\zeta|^{N }} e^{i\zeta\cdot z} d\zeta\Big| &\leq \frac {\|\eta'_0\|_{L^{\infty}}}{\sigma}\int_{B_1(x)}\frac{  d\zeta   }{|\zeta-x|^{N-1 }}
 \leq  \frac{2 \omega_{_{N-1}}}{\sigma}    ,
\end{align*}
since $\|\eta'_0\|_{L^\infty}\leq 2$.
This ends the proof. \hfill$\Box$\medskip

\noindent{\bf Proof of Theorem \ref{teo 1.1-sum}. } We recall that $\Phi_k(x,y)$ and $\widehat{\Phi}_k(z,y)$ have been defined in the proof of Theorem \ref{teo 1.3}. If we denote
$$\tilde v_{\sigma,z}(x):=v_{\sigma}(x,z)=w_\sigma(x)e^{{\rm i} x\cdot z},$$
 the projection of $v_{\sigma}$ onto the subspace of $L^2(\Omega)$ spanned by the $\phi_j$ for $1\leq j\leq k$ can be written in terms of the Fourier transform of $w_\sigma\Phi_k$ with respect to the $x$-variable:
$$\int_\Omega v_\sigma(x,z) \Phi_k(x,y) dx= (2\pi)^{N/2} \cF_x( w_\sigma\Phi_k)(z,y). $$
Put
$$v_{\sigma,k}(z,y)=v_{\sigma}(z,y)-(2\pi)^{N/2}\cF_x(w_\sigma \Phi_k)(z,y)  $$
and the Rayleigh-Ritz formula shows that
$$\lambda_{k+1}(\Omega) \int_{\Omega} |v_{\sigma,k}(z,y)|^2   dy \leq \int_{\Omega} \overline{v_{\sigma,k}(z,y)}  \lnlap_{,y} v_{\sigma,k}(z,y)dy $$
for any $z\in\R^N$ and $\sigma>0$, where the right hand side is a real value
$$\int_{\Omega} \overline{v_{\sigma,k}(z,y)}  \lnlap_{,y} v_{\sigma,k}(z,y)dy=\int_{\R^N} \overline{v_{\sigma,k}(z,y)}  \lnlap_{,y} v_{\sigma,k}(z,y)dy =\int_{\R^N} 2\ln|\xi| \Big| \cF(v_{\sigma,k})(z,\xi)\Big|^2 d\xi, $$
  although $v_{\sigma,k}$ is complex valued function.
Then, integrating this last inequality with respect to $z$ in $B_r\setminus B_1$, for $r>1$, 
we obtain
$$\displaystyle\lambda_{k+1}(\Omega)\leq \inf_{\sigma>0} \frac{\displaystyle\int_{B_r\setminus B_1} \int_{\Omega} \overline{v_{\sigma,k}(z,y)}  \lnlap_{,y} v_{\sigma,k}(z,y)dydz }{\displaystyle\int_{B_r\setminus B_1 } \int_{\Omega} |v_{\sigma,k}(z,y)|^2   dydz}.$$
By  Pythagore's theorem, we have that
$$\int_{\Omega} |v_{\sigma,k}(z,y)|^2   dy
 =\int_{\Omega} |v_{\sigma}(z,y)|^2   dy -(2\pi)^N\int_{\Omega}\sum^k_{j=1}|\cF_x(w_\sigma\phi_i)(z)|^2 \phi_i(y)^2 dy,
$$
  integrating over $B_r\setminus B_1$ implies that
\begin{align*}
&\int_{B_r\setminus B_1} \int_{\Omega} |v_{\sigma,k}(z,y)|^2   dydz\geq \frac{\omega_{_{N-1}}r^{N}}N \int_{\Omega} w_\sigma^2(y) dy-(2\pi)^N \sum^k_{j=1}\int_{B_r\setminus B_1} |\cF_x(w_\sigma\phi_i)(z)|^2dz.
 \end{align*}
On the other hand,
\begin{align*}
\int_{B_r\setminus B_1} \int_{\Omega} \overline{v_{\sigma,k}(z,y)}  \lnlap_{,y} v_{\sigma,k}(z,y)&   dydz  = \int_{B_r\setminus B_1} \int_{\Omega} \overline{v_{\sigma}(z,y)} \lnlap_{,y}  v_{\sigma}(z,y)    dydz
\\&  -(2\pi)^N
\int_{B_r\setminus B_1} \int_{\Omega}\overline{\cF_x(w_\sigma \Phi_k)(z,y)} \lnlap_{,y} \cF_x(w_\sigma \Phi_k)(z,y)  dy dz,
 \end{align*}
where
\begin{align*}
\int_{{B_r\setminus B_1}} \int_{\Omega}\overline{\cF_x(w_\sigma \Phi_k)(z,y) }\lnlap_{,y}  \cF_x(w_\sigma \Phi_k)(z,y)  dy dz = \sum^k_{j=1}\lambda_j(\Omega)\int_{{B_r\setminus B_1}} |\cF_x(w_\sigma \phi_j)(z) |^2dz
\end{align*}
and
\begin{align*}
 \int_{B_r\setminus B_1} \int_{\Omega} \overline{v_{\sigma}(z,y)} &\lnlap_{,y} v_{\sigma}(z,y)    dydz
\\ &\leq \int_{B_r\setminus B_1} \int_{\Omega} w_{\sigma}^2(y) |\lnlap_{,y}   e^{{\rm i} y.  z}|    dydz
 +\int_{B_r\setminus B_1} \int_{\Omega} w_{\sigma} (y)|\cL_z  w_{\sigma} (y) |    dydz
 \\&\leq\int_{B_r\setminus B_1} \int_{\Omega} w_{\sigma}^2(y)  \ln|z|   dydz+ \frac{2 \omega_{_{N-1}}}{\sigma} \int_{B_r\setminus B_1}\int_{\Omega}w_\sigma(y) dydz
 \\&= \frac{\omega_{_{N-1}}}N \varrho_{_{2,\sigma}}  \Big(r^N\ln r -\frac1N (r^N-1)\Big)
   +\frac{\omega_{_{N-1}}^2}{N \sigma}  \varrho_{_{1,\sigma}}  \Big( r^N-1 \Big)
   \\&\leq \frac{\omega_{_{N-1}}}N  \varrho_{_{2,\sigma}}  \,r^N\ln r + \frac{\omega_{_{N-1}}^2}{N \sigma}  \varrho_{1,\sigma}  \,   r^N
 \end{align*}
 with
 $$\varrho_{_{1,\sigma}}=\int_{\Omega} w_{\sigma}(y) dy\quad {\rm and}\quad \varrho_{_{2,\sigma}}=\int_{\Omega} w_{\sigma}^2(y) dy.$$
Because of Parseval's identity, there holds
 $$\int_{B_r\setminus B_1} |\cF_x(w_\sigma\phi_i)(z)|^2dz \leq \int_\Omega(w_\sigma\phi_i)^2 dx\leq 1. $$

Let $k_0$ be the smallest positive integer such that $k_0\geq \frac{e Nd_N}{2}|\Omega|$
and then $\lambda_{k_0}(\Omega)\geq 0$.

For  $k\geq k_0$, 
 we choose $r>1$ such that
$$\frac{ 2\omega_{_{N-1}}}{N^2} r^{N } \ln r \geq  \frac{\omega_{_{N-1}}r^{N}}N \Longleftrightarrow r\geq e^{\frac N2}\quad\text{
and} \quad \frac{\omega_{_{N-1}}r^{N}}N |\Omega|>2(2\pi)^N k,$$
 then we have that 
\begin{equation*} 
 \lambda_{k+1}(\Omega)\leq \frac{\frac{\omega_{_{N-1}} r^N}{N}\left(\varrho_{_{2,\sigma}} \frac 2N\ln\frac{r^N}{e}  + \varrho_{_{1,\sigma}} \frac{\omega_{_{N-1}}}{\sigma} \right)-(2\pi)^N\sum^k_{j=1}\lambda_j(\Omega)\int_{B_r} |\cF_x(w_\sigma \phi_j)(z) |^2dz}{\frac{\omega_{_{N-1}}r^{N}}N \varrho_{_{2,\sigma}} -(2\pi)^N \sum^k_{j=1}\int_{B_r} |\cF_x(w_\sigma\phi_j)(z)|^2dz}
 \end{equation*} 
Denote 
 $$A_1=\frac{\omega_{_{N-1}} r^N}{N}\left(\varrho_{_{2,\sigma}} \frac 2N\ln\frac{r^N}{e}  + \varrho_{_{1,\sigma}} \frac{\omega_{_{N-1}}}{\sigma} \right) \quad{\rm and}\quad A_2= \frac{\omega_{_{N-1}}r^{N}}N \varrho_{_{2,\sigma}},$$
 then 
\begin{align*}
0&\leq \frac{A_1-(2\pi)^N\sum^k_{j=1}\lambda_j(\Omega)\int_{B_r} |\cF_x(w_\sigma \phi_j)(z) |^2dz}{A_2 -(2\pi)^N \sum^k_{j=1}\int_{B_r} |\cF_x(w_\sigma\phi_j)(z)|^2dz}-\lambda_{k+1}(\Omega)
\\[2mm]&= \frac{\big(A_1-A_2\lambda_{k+1}(\Omega)\big)+(2\pi)^N\sum^k_{j=1}\Big(\lambda_{k+1}(\Omega)-\lambda_j(\Omega)\Big)\int_{B_r} |\cF_x(w_\sigma \phi_j)(z) |^2dz}{A_2 -(2\pi)^N \sum^k_{j=1}\int_{B_r} |\cF_x(w_\sigma\phi_j)(z)|^2dz}
\\[2mm]&\leq \frac{\big(A_1-A_2\lambda_{k+1}(\Omega)\big)+(2\pi)^N\sum^k_{j=1}\Big(\lambda_{k+1}(\Omega)-\lambda_j(\Omega)\Big) }{A_2 -(2\pi)^N k }
\end{align*} 
since $\lambda_{k+1}(\Omega)\geq \lambda_{j}(\Omega)$ for $j< k+1$ and $\int_{B_r} |\cF_x(w_\sigma\phi_j)(z)|^2dz\in(0,1)$. 
As a consequence, we obtain that 
 \begin{equation}\label{3.1-sum} 
\lambda_{k+1}(\Omega)\leq \frac{ \frac{\omega_{_{N-1}} r^N}{N}\left(\varrho_{_{2,\sigma}}\frac 2N\ln\frac{r^N}{e} + \varrho_{_{1,\sigma}}\frac{\omega_{_{N-1}}}{\sigma} \right)  -(2\pi)^N\sum^k_{j=1} \lambda_j(\Omega)   }{\frac{\omega_{_{N-1}}r^{N}}N \varrho_{_{2,\sigma}}-(2\pi)^N k}, 
 \end{equation} 
where
$$\frac{\omega_{_{N-1}}r^{N}}N \varrho_{_{2,\sigma}}-(2\pi)^N k>\frac{\omega_{_{N-1}}r^{N}}N \frac{|\Omega|}2-(2\pi)^N k>0.$$
  
\noindent We fix  $\sigma=\sigma_1$ and first impose $k,\, r>1$    such that
 $$\frac{\omega_{_{N-1}}r^{N}}N \varrho_{_{2,\sigma}} =(2\pi)^N (k+1),$$
and take $r=k^{\frac1N}$ for  $k\geq k_0 $, then  we recall that
$$\lambda_{k_0}\geq0 $$
and 
\begin{align*}
   (2\pi)^N\sum^{k+1}_{j=1} \lambda_j(\Omega)   
  & \leq   \frac{2\omega_{_{N-1}}}{N^2}  \varrho_{_{2,\sigma}}  \,r^N\ln \frac{r^N}e + \frac{\omega_{_{N-1}}^2}N \sigma^{-1}  |\Omega|^{\frac12} \sqrt{ \varrho_{_{2,\sigma}}} \,   r^N
\\&\leq   (2\pi)^N\frac{2(k+1)}{N}  \Big(\ln(k+1) +\ln \frac{p_N}{\varrho_{_{2,\sigma}}}+ \frac{\omega_{_{N-1}}}{\varrho_{_{2,\sigma}}}  \ln \ln \frac{p_N (k+1)}{\varrho_{_{2,\sigma}}} \Big)
\\&\leq  (2\pi)^N\frac{2(k+1)}{N}  \Big(\ln(k+1) +\ln \frac{2p_N}{|\Omega|}+ \frac{2\omega_{_{N-1}}}{\sqrt{|\Omega|}}  \ln \ln \frac{2p_N (k+1)}{|\Omega|} \Big),
\end{align*} 
where
$$p_N=\frac{2(2\pi)^NN}{\omega_{_{N-1}}} \quad{\rm and}\quad \frac{|\Omega|}2\leq \varrho_{_{2,\sigma}} \leq |\Omega|.$$

Moreover,  (\ref{upper bb-lim}) follows by the lower bound (\ref{lower bb}) and the upper bound (\ref{upper bb}) directly.  \hfill$\Box$\medskip

\section{ Further discussion }

\noindent{\bf Proof of Corollary  \ref{cr 1.2} $(i)$. } On the one hand, it follows by the nondecreasing monotonicity of $k\mapsto \lambda_k(\Omega)$ and (\ref{lower bb})  that
\begin{equation}\label{lower bb-1-y}
  \lambda_k(\Omega) \geq   \frac{1}{k}\sum^k_{i=1}\lambda_i(\Omega) \geq  \frac2{N} \left(\ln k+ \ln \Big(\frac{2}{eNd_N|\Omega|}\Big)-\ln\ln\Big(\frac{2k}{eNd_N|\Omega|}\Big)\right).
 \end{equation}

 On the other hand, we take $m=[\frac{k}{\ln\ln k}]+1$ and obtain
 \begin{align}
\lambda_{k+1}(\Omega) &  \leq \frac1m\Big(\sum^{k+m}_{j=1} \lambda_{j}(\Omega)-\sum^{k}_{j=1} \lambda_{j}(\Omega)\Big)\nonumber
\\& \leq   \frac{2}{Nm} (k+m)\left(\ln (k+m) +1+\ln \Big(\frac{p_N}{|\Omega|}\Big)+ \frac{\omega_{_{N-1}}}{\sqrt{|\Omega|}}  \ln \ln \Big(\frac{p_N (k+1+m)}{|\Omega|}\Big) \right)\nonumber
\\&\quad -  \frac{2}{Nm} k\left(\ln k+  \ln\Big(\frac{2}{eNd_N|\Omega|}\Big) -\ln\ln\Big(\frac{2k}{eNd_N|\Omega|}\Big)\right)\nonumber
\\& \leq   \frac{2}{N} \ln (k+m) + \frac{2}{N} (\ln\ln k)\ln\big(1+\frac1{\ln\ln k}\big)+ c_1\ln\ln k+c_2(\ln\ln k)^2\nonumber
\\&\leq  \frac{2}{N} \ln (k+1)+\frac{4}{N}  +c_3\big(\ln\ln (k+1)\big)^2,\label{upper bound -lower term}
\end{align}
where  $c_1,\, c_2,\, c_3>0$ is independent of $k$.
Thus, (\ref{lower bb-1}) holds true and we have that 
$$\lim_{k\to+\infty}\frac{\lambda_k(\Omega)}{\ln k}=\frac 2N.$$
We complete the proof. \hfill$\Box$

\subsection{ Proof of Corollary  \ref{cr 1.2} $(ii)$ }

With the help of some estimates in \cite{LW}, we can prove the bounds  (\ref{Weyl-upplow}) for $\lambda_k(\Omega)$. 

\begin{proposition}\label{pr 4.1-upp}
Let $\Omega$ be a bounded domain  in $\R^N$ and $ \lambda_k(\Omega) $  be the $k$-th eigenvalue of Dirichlet problem (\ref{eq 1.1}).   
  Then we have that for $k\geq    1$
$$\lambda_k(\Omega)\leq  \frac{2}{N} \ln k+c_3\big(\ln\ln (k+e)\big)^2+2\ln \frac1R, $$ 
where $R=\frac{|\Omega|}{|B_1|}$.
Particularly, 
$$\lambda_1(\Omega)\leq 2\ln \frac1R+c_3,$$
 where $c_3\in \R$ is independent of $\Omega$. 
\end{proposition}
\noindent{\bf Proof. }  Note that the constant $c_0$ in the second inequality (\ref{lower bb-1}) could is dependent of $|\Omega|$, but independent of the shape of $\Omega$. 
Let 
$$\Omega_R=\{R^{-\frac1N}x:\, x\in\Omega\}. $$ 
From (\ref{lower bb-1}) we can get that 
$$\lambda_{k}(\Omega_R)\leq  \frac{2}{N} \ln k+ c_4\big(\ln\ln (k+e)\big)^2,\quad k\geq \frac{e^{e+1}Nd_N}{2}|B_1|$$
where $c_4$ is independent of $k$ and $\Omega_R$.  The above inequality holds for $k\geq 1$ only by adjusting the constant $c_4$.  
It is shown a scaling property for logarithmic Laplacian \cite[Lemma 2.5]{LW} that 
$$\lambda_k(\Omega)=\lambda_k(t\Omega)-2\ln\frac1t,$$
where $t\Omega=\{tx:\, x\in\Omega\}$.
From the assumption that  $\Omega_R=R^{-1} \Omega$, we have that 
$$\lambda_k(\Omega) =\lambda_k(\Omega_R)+2\ln \frac1R. $$ 
We omit the left proof.\hfill$\Box$

\begin{proposition}\label{pr 4.2-low}
Let $\Omega$ be a bounded domain  in $\R^N$ and $ \lambda_k(\Omega) $  be the $k$-th eigenvalue of Dirichlet problem (\ref{eq 1.1}).   
  Then we have that for $k\geq   1$
$$\lambda_k(\Omega)\geq  \frac{2}{N} \Big(\ln k+\ln \frac{2}{eN d_N|\Omega|}\Big). $$ 
Particularly, $\lambda_1(\Omega)\geq \frac{2}{N} \Big(-\ln |\Omega|+c_5\Big)$, where $c_5=\ln \frac{2}{eN d_N} $. 
\end{proposition}
\noindent{\bf Proof. } It follows by\cite[Lemma 2.5]{LW} that 
$$\cN(\lambda)\leq \frac{|\Omega|\omega_{_{N-1}}}{N(2\pi)^N } e^{\frac{\lambda}{2}N+1}.$$
Note that 
$$\lambda_k(\Omega)\geq \lambda\ \ {\rm is\ equivalent\ to}\ \ \cN(\lambda)\leq k,$$
and then 
$$\lambda_k(\Omega)\geq \frac2N \ln \frac{N(2\pi)^N k } {|\Omega|\omega_{_{N-1}}}=\frac2N \ln \frac{2k}{eNd_N|\Omega|}.$$
Then all inequalities follow.  \hfill$\Box$


\bigskip

\bigskip

\medskip
 \noindent{\small {\bf Acknowledgements:}  Chen is is supported by NNSF of China, No: 12071189 and 12001252,
by the Jiangxi Provincial Natural Science Foundation, No: 20202BAB201005, 20202ACBL201001.

\end{document}